\documentclass[10pt]{amsart}
\usepackage{amssymb}
\usepackage[centertags]{amsmath}
\usepackage{amssymb}
\usepackage{xr}
\usepackage{amscd}
\usepackage{amsthm}
\usepackage{amsxtra}

\newtheorem{thm}{Theorem}[section]
\newtheorem{prop}[thm]{Proposition}
\newtheorem{lemma}[thm]{Lemma}
\newtheorem{rem}[thm]{Remark}
\newtheorem{cor}[thm]{Corollary}
\newtheorem{definition}[thm]{Definition}

\newtheorem{example}[thm]{Example}
\begin{document}
\author{Daniel Jupiter}
\address{Department of Mathematics\\ Texas A \& M University\\
College Station Texas, 77843-3368} \email{jupiter@math.tamu.edu}

\author{David Redett}
\address{Department of Mathematics\\ Texas A \& M University\\
College Station Texas, 77843-3368} \email{redett@math.tamu.edu}
\curraddr{Department of Mathematical Sciences\\ Indiana-Purdue
University Fort Wayne\\ Fort Wayne, IN 46805}
\email{redettd@ipfw.edu}


\date{November 2, 2004}


\begin{abstract}In this article we examine Dirichlet type spaces in the unit polydisc, and multipliers between these
spaces. These results extend the corresponding work of G. D. Taylor
in the unit disc.

In addition, we consider functions on the polydisc whose restrictions to lower dimensional polydiscs
lie in the corresponding Dirichet type spaces. We see that such functions need not be in the Dirichlet type
space of the whole polydisc. Similar observations are made regarding multipliers.

\end{abstract}

\title{Multipliers on Dirichlet type spaces}
\maketitle

\section{Introduction}
In \cite{taylor:multiplier} and \cite{taylor:thesis} G. D. Taylor
studied the Dirichlet type spaces, $D_{\alpha}$, on the unit disc in
$\mathbb{C}$. Additionally, he studied the multipliers between such
spaces: functions multiplying one space $D_{\alpha}$ to another
space $D_{\beta}$. An almost complete characterization of
multipliers was achieved.

In this article we generalize Taylor's results to the
case of Dirichlet type spaces on the unit polydisc in $\mathbb{C}^n$.
Following Taylor, we identify the power series of holomorphic functions in $D_{\alpha}$ with
elements in a weighted $\ell^2$ space. As our power series are multidimensional,
we allow the weights to be multidimensional as well. That is to say, the weight may vary from coordinate
direction to coordinate direction.

An interesting phenomenon arises when we restrict a function in
$D_{(\alpha_1,\,\ldots,\,\alpha_n)}$ on the unit polydisc to a
slice - a lower dimensional polydisc parallel to the coordinate
axes. Such a restricted function lies in the appropriate weighted
$\ell^2$ space of the lower dimensional polydisc. The converse to
this fact is, however, not true. In fact, a function may be in the
appropriate weighted $\ell^2$ space in each direction, without
being in the corresponding weighted $\ell^2$ space on the full
dimensional polydisc. This implies that Taylor's results cannot be
carried through to the higher dimensional case
merely by examining restrictions to slices.
These ideas are explored in Section \ref{section:slices}.\\

In the one dimensional setting the Dirichlet type spaces are
parameterized by $\alpha\in\mathbb{R}$. Much of Taylor's work
relies on carefully partitioning $\mathbb{R}$ and examining how
the spaces, and the multipliers between them, vary as the
parameter $\alpha$ changes, and moves from one element of the
partition of $\mathbb{R}$ to another.

As our weights may vary from coordinate to coordinate, our
parameter space is $\mathbb{R}^n$, rather than $\mathbb{R}$. We
can ``partition'' our parameter space in much the same way that
Taylor does. That is to say, we can insist that we look at those
parameter vectors where all the entries of the vector lie in one
element of Taylor's partition of $\mathbb{R}$. If we do so,
Taylor's results carry through. In fact, the same techniques of
proof can be employed in many instances, although there are cases
where alternate methods of proof need to be found.

The generalizations of Taylor's results to higher dimension can be
found in Section \ref{section:spaces}, which examines the spaces
$D_{\alpha}$, and in Section \ref{section:multipliers}, which
examines
multipliers between spaces.\\

For ease of notation we state and prove our results for the
bidisc. Most proofs carry through to higher dimensional polydiscs
with only the obvious changes necessary to account for differences
in dimension. In the few instances where this is not the case, we
comment on the changes needed for the proof in the higher
dimensional case.

\section{The Dirichlet Type Spaces, $D_{\alpha}$}\label{section:spaces}

\noindent Our objects of study are holomorphic functions on the
bidisc.

\begin{definition}[Holomorphic Function]
Let $\mathbb{U}^2$ denote the unit bidisc. We say that a function
\[f:\mathbb{U}^2\rightarrow \mathbb{C}\]
is holomorphic if it is holomorphic in each variable separately.
Each holomorphic function, $f$, on the bidisc can be represented
as a power series in $\mathbb{U}^2$:
\[f(z,\,w)=\sum_{(k,\,l)\in\mathbb{Z}^2_+}a_{k,\,l}z^kw^l,\] with $(z,\,w)\in\mathbb{U}^2$, and
$a_{k,\,l}\in\mathbb{C}$.
\end{definition}

\noindent We examine holomorphic functions satisfying a growth
condition on the coefficients of their Taylor series expansions.

\begin{definition}[$D_{\alpha}$]\label{definition:D}
Let $\alpha=(\alpha_{1},\,\alpha_{2})$ be an element of
$\mathbb{R}^2$. We define the space $D_{\alpha}$ as
\[D_{\alpha}=\Big\{f(z,\,w)=\sum_{(k,\,l) \in\mathbb{Z}^2_+}a_{k,\,l}z^kw^l \,;\,
\sum_{(k,\,l)\in\mathbb{Z}^2_+}
|a_{k,\,l}|^2(k+1)^{\alpha_{1}}(l+1)^{\alpha_{2}} <
\infty\Big\}.\]
\end{definition}

\begin{rem} As we shall see in Theorem \ref{thm:multipliers less than
0}, $D_{(0,\,0)}$ is the Hardy space $H^2(\mathbb{U}^2)$, and
$D_{(-1,\,-1)}$ is the Bergman space on $\mathbb{U}^2$.
\end{rem}

\noindent It is often useful to view $D_{\alpha}$ as a weighted
${\ell}^2$ space, with inner product defined as follows.

\begin{definition}[$\|f\|_{\alpha}$]
Let $f(z,\,w)=\sum_{(k,\,l)\in\mathbb{Z}^2_+}a_{k,\,l}z^kw^l$
and\\
$g(z,\,w)=\sum_{(k,\,l)\in\mathbb{Z}^2_+}b_{k,\,l}z^kw^l$ be two
elements of $D_{\alpha}$.  We define the inner product of $f$ and
$g$ as
\[(f,g)_{\alpha}=\sum_{(k,\,l)\in\mathbb{Z}^2_+}
a_{k,\,l}\overline{b_{k,\,l}}(k+1)^{\alpha_{1}}(l+1)^{\alpha_{2}}.\]

The norm of an element, $f$, in $D_{\alpha}$ is denoted by
$\|f\|_{\alpha}$ and given as
\[\|f\|_{\alpha}=(f,\,f)_{\alpha}^{\frac{1}{2}}=\biggl(\sum_{(k,\,l)\in\mathbb{Z}^2_+}
|f_{k,\,l}|^2(k+1)^{\alpha_{1}}(l+1)^{\alpha_{2}}\biggl)^{\frac{1}{2}}.\]
\end{definition}

\noindent Definition \ref{definition:D} does not stipulate that
elements of $D_{\alpha}$ are holomorphic, and it is not, in fact,
a priori clear that the functions in $D_{\alpha}$ are holomorphic.
We have, however, the following proposition.

\begin{prop}\label{prop:D is holo}
Let $\alpha=(\alpha_1,\,\alpha_2)\in\mathbb{R}^2$. If
$\{a_{k,\,l}\}_{(k,\,l)\in\mathbb{Z}^2_+}$ is a sequence in
$\mathbb{C}$ satisfying
\[\sum_{(k,\,l)\in\mathbb{Z}^2_+}|a_{k,\,l}|^2(k+1)^{\alpha_{1}}(l+1)^{\alpha_{2}}
< \infty,\] then
$f(z,\,w)=\sum_{(k,\,l)\in\mathbb{Z}^2_+}a_{k,\,l}z^kw^l$ converges
uniformly on compact subsets in $\mathbb{U}^2$, and is thus
holomorphic on $\mathbb{U}^2$.
\end{prop}

\noindent The proof of Proposition \ref{prop:D is holo} relies on
Abel's Lemma.

\begin{lemma}[Abel's Lemma]
Let $\{a_{k,\,l}\}_{(k,\,l)\in\mathbb{Z}^2_+}$ be a sequence in
$\mathbb{C}$, and let $(z_{0},\,w_{0})$ be a point in
$\mathbb{C}^2$. Assume that
\[\sup_{(k,\,l)\in\mathbb{Z}^2_+}|a_{k,\,l}z_{0}^kw_{0}^l| < M\] for some $M < \infty$.  Then
\[\sum_{(k,\,l)\in\mathbb{Z}^2_+}a_{k,\,l}z^kw^l\] converges uniformly
on compact subsets in $|z_{0}|\mathbb{U} \times |w_{0}|\mathbb{U}$.
\end{lemma}

\begin{proof}[Proof of Proposition \ref{prop:D is holo}]
We note first that the convergence of
\[\sum_{(k,\,l)\in\mathbb{Z}^2_+}|a_{k,\,l}|^2(k+1)^{\alpha_{1}}(l+1)^{\alpha_{2}}\]
implies that the terms
$|a_{k,\,l}|^2(k+1)^{\alpha_{1}}(l+1)^{\alpha_{2}}$ are bounded.

Next, we fix $0\leq r_{1},\,r_{2}<1$ and show that
\[\sup_{(k,\,l)\in\mathbb{Z}^2_+}|a_{k,\,l}|r_{1}^kr_{2}^l<M<\infty.\]
Then, by Abel's Lemma, we have that
$f(z,\,w)=\sum_{(k,\,l)\in\mathbb{Z}^2_+}a_{k,\,l}z^kw^l$ converges
uniformly on compact subsets in the bidisc $r_{1}\mathbb{U} \times
r_{2}\mathbb{U}$.  Since $r_1$ and $r_2$ were arbitrary, we conclude
that $f(z,\,w)$ converges uniformly on compact subsets in
$\mathbb{U}^2$.

There exists $K$ such that for $k,\,l > K$ we have
\begin{eqnarray}r_{1}^{k}r_{2}^l <
(k+1)^{\alpha_{1}/2}(l+1)^{\alpha_{2}/2}.\label{inequality:D is
holo}
\end{eqnarray}

Thus
\[|a_{k,\,l}|r_{1}^{k}r_{2}^l <
|a_{k,\,l}|(k+1)^{\alpha_{1}/2}(l+1)^{\alpha_{2}/2},\] for all
$k,\,l> K$.

From the remark at the beginning of this proof we see that the
numbers $|a_{k,\,l}|(k+1)^{\alpha_{1}/2}(l+1)^{\alpha_{2}/2}$ are
bounded. Hence, $|a_{k,\,l}|r_{1}^{k}r_{2}^l$ is bounded for
$k,\,l>K$.

There exists $N_1$ such that if $k> N_1$ and $0\leq l\leq K$, then
(\ref{inequality:D is holo}) holds. Reasoning as above, we
conclude that $|a_{k,\,l}|r_{1}^{k}r_{2}^l$ is bounded for $0 \leq
l \leq K$ and $k > N_1$.  An analogous argument gives us the
existence of $N_{2}$ such that $|a_{k,\,l}|r_{1}^{k}r_{2}^l$ is
bounded for $l > N_{2}$ and $0 \leq k \leq K$.

Only finitely many terms have not been estimated:
$|a_{k,\,l}|r_{1}^{k}r_{2}^l$, with $0\leq k\leq N_1$ and $0\leq l
\leq N_2$. These are, of course, bounded. We conclude that there
is an $M<\infty$ such that
$\sup_{(k,\,l)\in\mathbb{Z}^2_+}|a_{k,\,l}|r_{1}^kr_{2}^l<M.$
\end{proof}

\noindent We next define a useful linear functional on
$D_{\alpha}$: point evaluation.

\begin{definition}[$\lambda_{(z,\,w)}^{\alpha}$]
Let $f$ be an element of $D_{\alpha}$. We define
$\lambda_{(z,\,w)}^{\alpha}(f)$ as
\[\lambda_{(z,\,w)}^{\alpha}(f)=f(z,\,w).\]
\end{definition}

\noindent While it is clear that $\lambda_{(z,\,w)}^{\alpha}$ is a
linear functional on $D_{\alpha}$, we will need to use the
following technical lemma to show that it is bounded.

\begin{lemma}\label{lemma:multiindex convergence}
For all $(z,\,w)$ in $\mathbb{U}^2$, the sum
\[\sum_{(k,\,l)\in\mathbb{Z}^2_+}|z|^{2k}|w|^{2l}(k+1)^{-\alpha_{1}}(l+1)^{-\alpha_{2}}\]
converges.
\end{lemma}

\begin{proof}
Fix $(z,\,w)$ in $\mathbb{U}^2$.  There exists $K$ such that for
$k,\,l > K$ we have

\[|z|^{2k}|w|^{2l} <(k+1)^{\alpha_{1}-2}(l+1)^{\alpha_{2}-2}.\]

Thus we have
\begin{eqnarray}\label{inequality:multiindex convergence}
|z|^{2k}|w|^{2l}(k+1)^{-\alpha_{1}}(l+1)^{-\alpha_{2}}
<(k+1)^{-2}(l+1)^{-2}\end{eqnarray} for $k,\,l> K$.

For $0 \leq l \leq K$, there exists $N_{1}$ such that if $k>
N_{1}$, then (\ref{inequality:multiindex convergence}) holds.
Analogously, for $0 \leq k \leq K$, there exists $N_{2}$ so that
if $l > N_{2}$, then (\ref{inequality:multiindex convergence})
holds.

Let $N=\max\{K,\,N_{1},\,N_{2}\}$. Then
\begin{align*}
\sum_{(k,\,l)\in\mathbb{Z}^2_+}|z|^{2k}|w|^{2l}&(k+1)^{-\alpha_{1}}(l+1)^{-\alpha_{2}}\\
&\leq \sum_{0 \leq l,\,k \leq N}|z|^{2k}|w|^{2l}(k+1)^{-\alpha_{1}}(l+1)^{-\alpha_{2}}\\
&+\sum_{(k,\,l)\in\mathbb{Z}^2_+}(k+1)^{-2}(l+1)^{-2}.
\end{align*}
Since
\[\sum_{(k,\,l)\in\mathbb{Z}^2_+}(k+1)^{-2}(l+1)^{-2}\]
converges, we have that
\[\sum_{(k,\,l)\in\mathbb{Z}^2_+}|z|^{2k}|w|^{2l}(k+1)^{-\alpha_{1}}(l+1)^{-\alpha_{2}}\]
converges.
\end{proof}

\begin{lemma}\label{lemma:point evaluation}
The functional $\lambda_{(z,\,w)}^{\alpha}$ is a bounded linear
functional on $D_{\alpha}$, with norm
\[\Big(\sum_{(k,\,l)\in\mathbb{Z}^2_+}|z|^{2k}|w|^{2l}(k+1)^{-\alpha_{1}}(l+1)^{-\alpha_{2}}\Big)^{1/2}.\]
\end{lemma}

\begin{proof}
Let $f(z,\,w)=\sum_{(k,\,l)\in\mathbb{Z}^2_+}a_{k,l}z^kw^l\in
D_{\alpha}$. Then
\begin{align*}
|\lambda_{(z,w)}^{\alpha}(f)| & = |f(z,w)| \\
&\leq \sum_{(k,\,l)\in\mathbb{Z}^2_+}|a_{k,\,l}||z|^k|w|^l \\
&=\sum_{(k,\,l)\in\mathbb{Z}^2_+}|a_{k,\,l}||z|^k|w|^l(k+1)^{\alpha_{1}/2}(l+1)^{\alpha_{2}/2}
      (k+1)^{-\alpha_{1}/2}(l+1)^{-\alpha_{2}/2} \\
&\leq
\Big(\sum_{(k,\,l)\in\mathbb{Z}^2_+}|a_{k,\,l}|^2(k+1)^{\alpha_{1}}(l+1)^{\alpha_{2}}\Big)^{1/2}\\
 &\quad\quad \cdot
 \Big(\sum_{(k,\,l)\in\mathbb{Z}^2_+}|z|^{2k}|w|^{2l}(k+1)^{-\alpha_{1}}(l+1)^{-\alpha_{2}}\Big)^{1/2}\\
&=
\|f\|_{\alpha}\Big(\sum_{(k,\,l)\in\mathbb{Z}^2_+}|z|^{2k}
|w|^{2l}(k+1)^{-\alpha_{1}}(l+1)^{-\alpha_{2}}\Big)^{1/2}.
\end{align*}

We see that $\|\lambda_{(z,w)}^{\alpha}\| \leq
\Big(\sum_{(k,\,l)\in\mathbb{Z}^2_+}|z|^{2k}|w|^{2l}(k+1)^{-\alpha_{1}}
(l+1)^{-\alpha_{2}}\Big)^{1/2}$.\\

Let $(z_{0},w_{0})$ be in $\mathbb{U}^2$, and define
$K_{(z_0,\,w_0)}^{\alpha}:\mathbb{U}^2\rightarrow\mathbb{C}$ by
\[K_{(z_{0},\,w_{0})}^{\alpha}(z,\,w) = \sum_{(k,\,l)\in\mathbb{Z}^2_+}
\overline{z_{0}}^k\overline{w_{0}}^lz^kw^l(k+1)^{-\alpha_{1}}(l+1)^{-\alpha_{2}}.\]
We note that $K_{(z_0,\,w_0)}^{\alpha}$ is in $D_{\alpha}$ since
\[\|K_{(z_0,\,w_0)}^{\alpha}\|_{\alpha}^2= \sum_{(k,\,l)\in\mathbb{Z}^2_+}|z_0|^{2k}|w_0|^{2l}
(k+1)^{-\alpha_{1}}(l+1)^{-\alpha_{2}}\] converges, by Lemma
\ref{lemma:multiindex convergence}.

Now,

\begin{align*}
|\lambda_{(z,\,w)}^{\alpha}(K_{(z,\,w)}^{\alpha})| & =
\sum_{(k,\,l)\in\mathbb{Z}^2_+}|z|^{2k}|w|^{2l}(k+1)^{-\alpha_{1}}(l+1)^{-\alpha_{2}} \\
=& \Big(
\sum_{(k,\,l)\in\mathbb{Z}^2_+}|z|^{2k}|w|^{2l}(k+1)^{-\alpha_{1}}(l+1)^{-\alpha_{2}}\Big)^{1/2}\\
&\quad\quad\cdot\Big( \sum_{(k,\,l)\in\mathbb{Z}^2_+}|z|^{2k}|w|^{2l}(k+1)^{-\alpha_{1}}(l+1)^{-\alpha_{2}}\Big)^{1/2} \\
&= \Big(
\sum_{(k,\,l)\in\mathbb{Z}^2_+}|z|^{2k}|w|^{2l}(k+1)^{-\alpha_{1}}(l+1)^{-\alpha_{2}}\Big)^{1/2}
\|K_{(z,w)}^{\alpha}\|_{\alpha}.
\end{align*}

Hence $\|\lambda_{(z,\,w)}^{\alpha}\| = \Big(
\sum_{(k,\,l)\in\mathbb{Z}^2_+}|z|^{2k}|w|^{2l}(k+1)^{-\alpha_{1}}(l+1)^{-\alpha_{2}}\Big)^{1/2}$.
\end{proof}

\begin{rem}
In the above proof, $K_{(z,\,w)}^{\alpha}$ is the reproducing
kernel for $D_{\alpha}$.
\end{rem}

\noindent In order to characterize and compare Dirichlet type
spaces, we define a partial ordering on our multiindices.
\begin{definition}[$\alpha \succ \beta$, $\alpha\succeq\beta$]
Let $\alpha=(\alpha_1,\,\alpha_2)$ and $\beta=(\beta_1,\,\beta_2)$.
We write $\alpha \succ \beta$  if $\alpha_{1}
> \beta_{1}$ and $\alpha_{2} > \beta_{2}$, and $\alpha \succeq \beta$  if $\alpha_{1}
\geq \beta_{1}$ and $\alpha_{2} \geq \beta_{2}$.
\end{definition}

\noindent We make some basic observations about Dirichlet type
spaces, and the relationships between them.

\begin{prop}\label{prop:D in H infty}
If $\alpha \succ (1,\,1)$, then $D_{\alpha} \subseteq
H^\infty(\mathbb{U}^2)$: the set of bounded holomorphic functions
on $\mathbb{U}^2$.
\end{prop}

\begin{proof}
We have that
 \[\sum_{(k,\,l)\in\mathbb{Z}^2_+}|z|^{2k}|w|^{2l}
(k+1)^{-\alpha_{1}}(l+1)^{-\alpha_{2}} \leq
\sum_{(k,\,l)\in\mathbb{Z}^2_+}(k+1)^{-\alpha_{1}}(l+1)^{-\alpha_{2}}.\]
This last sum is finite if $\alpha_{1},\,\alpha_2>1$. This
observation, together with Lemma \ref{lemma:point evaluation},
finishes the proof.
\end{proof}

\noindent We remark that the containment here must be strict, as
$D_{\alpha}$ is a Hilbert space while $H^{\infty}(\mathbb{U})$ is
not.

\noindent It is clear that if $\alpha\succeq \beta$, then
$D_{\alpha}\subseteq D_{\beta}$. If $\alpha\succ\beta$ we have
strict containment.

\begin{prop}\label{prop:proper containment}
If $\beta \prec \alpha$ then $D_{\alpha} \subset D_{\beta}$.
\end{prop}

\begin{proof}
Let
\[f(z,\,w)=\sum_{(k,\,l)\in\mathbb{Z}^2_+}(k+1)^{(-\alpha_{1}-1)/2}
(l+1)^{(-\alpha_{2}-1)/2}z^kw^l.\] We see that
\[\|f\|_{\beta}^2=\sum_{(k,\,l)\in\mathbb{Z}^2_+}(k+1)^{-\alpha_{1}+
\beta_{1}-1}(l+1)^{-\alpha_{2}+\beta_{2}-1}.\] Since
$-\alpha_{j}+\beta_{j}-1<-1$ for $j=1,\,2$, this last sum is
finite, and we conclude that $f$ is in $D_{\beta}$.

On the other hand
\[\|f\|_{\alpha}^2=\sum_{(k,\,l)\in\mathbb{Z}^2_+}(k+1)^{-1}(l+1)^{-1} =
\infty.\] Thus $f$ is not in $D_{\alpha}$, and we conclude that
$D_{\alpha}\subset D_{\beta}$.
\end{proof}

\begin{rem}\label{rem:proper containment} We can, in fact, say more. If
$\alpha=(\alpha_1,\,\alpha_2)\succeq\beta=(\beta_1,\,\beta_2)$,
with $\alpha_2=\beta_2$ but $\alpha_1>\beta_1$, then
$D_{\alpha}\subset D_{\beta}$. To see this, consider the function
\[f(z,\,w)=\sum_{k\in\mathbb{Z}_+}(k+1)^{(-\alpha_{1}-1)/2}
z^k.\] Just as in Proposition \ref{prop:proper containment}, we
see that $f$ is in $D_{\beta}$ but not in $D_{\alpha}$.
\end{rem}

\noindent We conclude this section by noting that if $\alpha\prec
(0,\,0)$ or $\alpha=(0,\,0)$ then the norm $\|f\|_{\alpha}$ is
equivalent to an integral. The space $D_{(0,\,0)}$ is the Hardy
space $H^2(\mathbb{U}^2)$, and thus has norm
\begin{equation}\label{equation:hardy}
\|f\|_2=\biggl(\frac{1}{(2\pi)^2}\sup_{0\leq
r<1}\int_{0}^{2\pi}\int_{0}^{2\pi}|f(re^{i\theta_1},\,re^{i\theta_2})|^2\,d\theta_1
d\theta_2\biggr)^{\frac{1}{2}}.
\end{equation}
(See, e.g. \cite{rudin:polydisc}.) This is, in fact, exactly the
same as $\|f\|_{(0,\,0)}$. If $\alpha\prec(0,\,0)$ we have the
following.

\begin{lemma}\label{lemma:integral norm equivalence}
For $\alpha \prec (0,0)$ the norm, $\|f\|_{\alpha}$, of
$D_{\alpha}$ is equivalent to the integral
\begin{align*}
\biggl(\frac{1}{\pi^2}\int_{0}^{1}&\int_{0}^{1}\int_{0}^{2\pi}\int_{0}^{2\pi}|
f(r_{1}e^{i\theta_{1}},\, r_{2}e^{i\theta_{2}})|^2(1-r_{1}^2)^{-1-\alpha_{1}}(1-r_{2}^2)^{-1-\alpha_{2}}\\
&\quad\quad r_1r_2\,d\theta_{1}d\theta_{2}dr_1dr_2.\biggr)^{\frac{1}{2}}
\end{align*}
\end{lemma}

\begin{proof}
We examine the above integral.

\begin{align*}
\frac{1}{\pi^2}\int_{0}^{1}\int_{0}^{1}&\int_{0}^{2\pi}\int_{0}^{2\pi}|f(r_{1}e^{i\theta_{1}},
r_{2}e^{i\theta_{2}})|^2(1-r_{1}^2)^{-1-\alpha_{1}}(1-r_{2}^2)^{-1-\alpha_{2}}\\
&\hspace{2in} r_1r_{2}\,d\theta_{1}d\theta_{2}dr_1dr_2 \\
 &\hspace{.3in}=4\int_{0}^{1}\int_{0}^{1}\sum_{(k,\,l)\in\mathbb{Z}^2_+}
|a_{k,\,l}|^2(1-r_{1}^2)^{-1-\alpha_{1}}(1-r_{2}^2)^{-1-\alpha_{2}}\\
&\hspace{2in} r_{1}^{2k+1}r_{2}^{2l+1}\,dr_{1}dr_{2} \\
&\hspace{.3in}= 4\sum_{(k,\,l)\in\mathbb{Z}^2_+}|a_{k,\,l}|^2\int_{0}^{1}\int_{0}^{1}(1-r_{1}^2)^{-1-\alpha_{1}}
(1-r_{2}^2)^{-1-\alpha_{2}}\\
&\hspace{2in} r_{1}^{2k+1}r_{2}^{2l+1}\,dr_{1}dr_{2} \\
&\hspace{.3in}= 4\sum_{(k,\,l)\in\mathbb{Z}^2_+}|a_{k,\,l}|^2(-\alpha_{1})^{-1}(^{k-\alpha_{1}}_{-\alpha_{1}})^{-1}
(-\alpha_{2})^{-1}(^{l-\alpha_{2}}_{-\alpha_{2}})^{-1}.
\end{align*}

The last equality follows from an integration by parts.

Hardy \cite[Chapter 5]{hardy:divergent} shows that
\[(-\alpha_{1})^{-1}(^{k-\alpha_{1}}_{-\alpha_{1}})^{-1}\]
is asymptotic to
\[(k+1)^{\alpha_{1}}.\]

Using this fact  and an argument similar to that in the proof of
Proposition \ref{prop:D is holo}, we conclude that
\[4\sum_{(k,\,l)\in\mathbb{Z}^2_+}|a_{k,l}|^2(-\alpha_{1})^{-1}(^{k-\alpha_{1}}_{-\alpha_{1}})^{-1}(-\alpha_{2})^{-1}
(^{l-\alpha_{2}}_{-\alpha_{2}})^{-1}\]
is comparable to
\[\sum_{(k,\,l)\in\mathbb{Z}^2_+}|a_{k,\,l}|^2(k+1)^{\alpha_{1}}(l+1)^{\alpha_{2}}.\]
\end{proof}

\begin{rem}
We note that in the above proof we can in fact allow $\alpha_1=0$
and $\alpha_2<0$. To see this, we replace the integral
\[\frac{1}{\pi}\int_0^{2\pi}\int_0^1|f(r_{1}e^{i\theta_{1}},
r_{2}e^{i\theta_{2}})|^2(1-r_{1}^2)^{-1-\alpha_{1}} \,dr_1d\theta_1\]
with
the supremum
\[\frac{1}{2\pi}\sup_{0\leq r_1<1}\int_0^{2\pi}|f(r_{1}e^{i\theta_{1}},
r_{2}e^{i\theta_{2}})|^2d\theta_1.\]
The proof then proceeds in essentially the same manner as above.

We conclude
that for $\alpha\preceq(0,\,0)$, the norm of $D_{\alpha}$ is
equivalent to an integral.
\end{rem}

\section{Multipliers}\label{section:multipliers}
\noindent A multiplier from $D_{\alpha}$ to $D_{\beta}$ is a
function which acts by pointwise multiplication as a mapping from
$D_{\alpha}$ to $D_{\beta}$. We remark that multipliers often
arise in functional analytic contexts; for example in the theory
of integral and differential operators.

\noindent We provide a precise definition.

\begin{definition}[Multiplier, $\mathcal{M}(D_{\alpha},\,D_{\beta})$]
Let $\mathcal{M}(D_{\alpha},\,D_{\beta})$ be defined as
\[\mathcal{M}(D_{\alpha},\,D_{\beta})=\{\phi:\mathbb{U}^2\rightarrow
\mathbb{C}\,;\, \phi f \in D_{\beta}, \textrm{ for all } f \in
D_{\alpha}\}.\]

We write $\mathcal{M}(D_{\alpha})$ for
$\mathcal{M}(D_{\alpha},\,D_{\alpha})$

An element of $\mathcal{M}(D_{\alpha},\,D_{\beta})$ (resp.
$\mathcal{M}(D_{\alpha})$) is called a multiplier from
$D_{\alpha}$ to $D_{\beta}$ (resp. a multiplier of $D_{\alpha}$).
\end{definition}

\noindent Our goal is to obtain as complete a characterization as possible of the
spaces $\mathcal{M}(D_{\alpha},\,D_{\beta})$ and $\mathcal{M}(D_{\alpha})$.
We begin with some basic remarks about multipliers.\\

\noindent We note first that since $1$ is in $D_{\alpha}$ for all
$\alpha$ in $\mathbb{R}^2$, we have that
$\mathcal{M}(D_{\alpha},D_{\beta}) \subseteq D_{\beta}$. Next, for
$h$ in $\mathcal{M}(D_{\alpha},D_{\beta})$, we define $T_{h} :
D_{\alpha} \rightarrow D_{\beta}$ by $$T_{h}f = hf,$$ for $f$ in
$D_{\alpha}$.  It is clear that $T_{h}$ is a linear
transformation; in fact, $T_h$ is bounded.

\begin{lemma}\label{lemma:T bounded}
For $h$ in $\mathcal{M}(D_{\alpha},D_{\beta})$, $T_{h}$ is a
bounded linear transformation.
\end{lemma}

\noindent Lemma \ref{lemma:point evaluation} states that point
evaluation is a bounded linear functional on $D_{\alpha}$. It is
well known that a multiplier between two Banach spaces where point
evaluation is bounded is a bounded linear transformation. See,
e.g. \cite[Chapter 1, Section 1]{taylor:thesis}.

\noindent We denote the norm of $T_h$ by
$\|T_h\|_{\alpha,\,\beta}$, suppressing the $\alpha$ and $\beta$,
if they are clear from the context. The lemma tells us that given
$h\in\mathcal{M}(D_{\alpha},\,D_{\beta})$, we obtain a bounded
linear transformation. We have, in addition, pointwise estimates
on elements of $\mathcal{M}(D_{\alpha},\,D_{\beta})$.

\begin{lemma}\label{lemma:multiplier bound}
Let $h$ be in $\mathcal{M}(D_{\alpha},\,D_{\beta})$. Then
\[|h(z,\,w)| \leq \|T_{h}\|_{\alpha,\,\beta}\frac{\|\lambda_{(z,\,w)}^{\beta}\|}{\|\lambda_{(z,\,w)}^{\alpha}\|}.\]
\end{lemma}

\noindent This lemma is also a well known result, following from the boundedness of point
evaluation. See again \cite[Chapter 1, Section 1]{taylor:thesis}.

\noindent As a simple consequence of this lemma, we obtain that if $\alpha\preceq\beta$ then the multipliers
from $D_{\alpha}$ to $D_{\beta}$ are bounded analytic functions.

\begin{cor}\label{corollary:H infty multipliers}
If $\alpha \preceq \beta$ then $M(D_{\alpha},\,D_{\beta}) \subseteq H^\infty(\mathbb{U}^2)$.
In particular, $M(D_{\alpha}) \subseteq H^\infty(\mathbb{U}^2)$.
\end{cor}

\begin{proof}
If $\alpha \preceq \beta$ then $\|\lambda_{(z,\,w)}^{\beta}\| \leq \|\lambda_{(z,\,w)}^{\alpha}\|$. Applying
Lemma \ref{lemma:multiplier bound} yields the result.
\end{proof}

\begin{rem} We shall see later (Theorem \ref{thm:multipliers alpha <
beta} and Proposition \ref{prop:beta > alpha}) that
$M(D_{\alpha},\,D_{\beta})=\{0\}$ if $\alpha \prec\beta$.
\end{rem}

\noindent The following interpolation result will be central in some of our later estimates.

\begin{thm}\label{theorem:interpolation}
 Let
\begin{align*}
\alpha^1=(\alpha^1_1,\,\alpha^1_2) &\textrm{ and } \beta^1=(\beta^1_1,\,\beta^1_2), \textrm{ and }\\
\alpha^2=(\alpha^2_1,\,\alpha^2_2) &\textrm{ and } \beta^2=(\beta^2_1,\,\beta^2_2).\\
\end{align*}
Let
\[\alpha=(\alpha_1,\,\alpha_2)=((1-\lambda)\alpha^1_1+\lambda\alpha^2_1,\,(1-\lambda)\alpha^1_2+\lambda\alpha^2_2)\]
and
\[\beta=(\beta_1,\,\beta_2)=((1-\lambda)\beta^1_1+\lambda\beta^2_1,\,(1-\lambda)\beta^1_2+\lambda\beta^2_2).\]
If
\[h\in\mathcal{M}(D_{\alpha^1},\,D_{\beta^1})\cap
\mathcal{M}(D_{\alpha^2},\,D_{\beta^2}),\] then
$h\in\mathcal{M}(D_{\alpha},\,D_{\beta})$ and
\[\|T_h\|_{\alpha,\,\beta}\leq\|T_h\|_{\alpha^1,\,\beta^1}^{1-\lambda}\|T_h\|_{\alpha^2,\,\beta^2}^{\lambda}.\]
\end{thm}

\noindent The proof of Theorem \ref{theorem:interpolation} is
essentially the same as Taylor's proof of the analogous result in
the one dimensional case, except that we are interpolating between
points in $\mathbb{R}^2$ rather than between points in $\mathbb{R}$.
The proofs are, in fact, so similar that we refer the reader to
Taylor's proof \cite[Chapter 1, Section 3]{taylor:thesis}.  We point
out, as Taylor did in \cite{taylor:multiplier} that the proof is
very similar to the proof of the Riesz-Thorin Interpolation Theorem.

Further, following a suggestion from the referee, we take a moment
to comment on the ``geometric interpolation spaces''. The idea is as
follows. Let
$(\mathcal{K}_{0},\,\mathcal{L}_0)=(D_{\alpha^1},\,D_{\beta^1})$,
$(\mathcal{K}_{1},\,\mathcal{L}_1)=(D_{\alpha^2},\,D_{\beta^2})$,
and $T$ be a bounded linear transformation that maps
$\mathcal{K}_{0}$ to $\mathcal{L}_0$ and $\mathcal{K}_{1}$ to
$\mathcal{L}_1$. We would like to find a scale of Hilbert spaces
$(\mathcal{K}_{\lambda},\,\mathcal{L}_{\lambda})$, $(0<\lambda<1)$,
connecting $(\mathcal{K}_{0},\,\mathcal{L}_0)$ and
$(\mathcal{K}_{1},\,\mathcal{L}_1)$, such that $T$ maps
$\mathcal{K}_{\lambda}$ to $\mathcal{L}_{\lambda}$ and
\[\|T\|_{\mathcal{K}_{\lambda}\rightarrow\mathcal{L}_{\lambda}}\leq
\|T\|_{\mathcal{K}_{0}\rightarrow\mathcal{L}_{0}}^{1-\lambda}\|T\|_{\mathcal{K}_{1}\rightarrow\mathcal{L}_{1}}^{\lambda}.\]
The spaces $(\mathcal{K}_{\lambda},\,\mathcal{L}_{\lambda})$ are
called the geometric interpolation spaces. In \cite[Appendix C]{A-M}
it is shown that if $\mathcal{K}_0$ and $\mathcal{K}_1$ are
``compatible'', as are $\mathcal{L}_0$ and $\mathcal{L}_1$, then the
geometric interpolation spaces exist, are unique, and a method to
construct them is given. Following this method, one sees that in our
case the geometric interpolation spaces are precisely the
$(D_{\alpha},\,D_{\beta})$, where
\[\alpha=(\alpha_1,\,\alpha_2)=((1-\lambda)\alpha^1_1+\lambda\alpha_1^2,\,(1-\lambda)\alpha_2^1
+\lambda\alpha_2^2)\] and
\[\beta=(\beta_1,\,\beta_2)=((1-\lambda)\beta^1_1+\lambda\beta_1^2,\,(1-\lambda)\beta_2^1
+\lambda\beta_2^2),\] with $0<\lambda<1$.

\begin{rem} We also easily see that if $h\in\mathcal{M}(D_{\alpha},\,D_{\beta})$, then
$h\in\mathcal{M}(D_{\alpha},\,D_{\mu})$ for $\mu\preceq\beta$. This follows from the fact
that $D_{\beta}\subseteq D_{\mu}$ if $\mu\preceq\beta$.

Similarly, if $h\in\mathcal{M}(D_{\alpha},\,D_{\beta})$, then
$h\in\mathcal{M}(D_{\tau},\,D_{\beta})$ for $\tau\succeq\alpha$.
\end{rem}

\noindent We are now in a position to characterize some of the spaces
$\mathcal{M}(D_{\alpha},\,D_{\beta})$.

\begin{thm}\label{thm:multipliers less than 0}
Let $(0,\,0)\succ\alpha\succeq\beta$. Then
$h\in\mathcal{M}(D_{\alpha},\,D_{\beta})$ if and only if
\[|h(z,\,w)|\leq
C\biggl((1-|z|^2)^{\frac{\beta_1-\alpha_1}{2}}(1-|w|^2)^{\frac{\beta_2-\alpha_2}{2}}\biggr).\]
\end{thm}

\begin{proof}
Suppose that $h$ is in $M(D_{\alpha},D_{\beta})$. Then
\begin{align*}
|h(z,\,w)| &\leq
\|T_{h}\|\frac{\|\lambda_{(z,\,w)}^{\beta}\|}{\|\lambda_{(z,\,w)}^{\alpha}\|} \\
&=\|T_{h}\|\frac{\Big(\sum_{(k,\,l)\in\mathbb{Z}^2_+}|z|^{2k}|w|^{2l}(k+1)^{-\beta_{1}}
(l+1)^{-\beta_{2}}\Big)^{1/2}}{\Big(\sum_{(k,\,l)\in\mathbb{Z}^2_+}|z|^{2k}|w|^{2l}
(k+1)^{-\alpha_{1}}(l+1)^{-\alpha_{2}}\Big)^{1/2}} \\
&=\|T_{h}\|\biggl(\frac{\sum_{k\in\mathbb{Z}_+}|z|^{2k}(k+1)^{-\beta_{1}}
\sum_{l\in\mathbb{Z}_+} |w|^{2l}
(l+1)^{-\beta_{2}}}{\sum_{k\in\mathbb{Z}_+}|z|^{2k}(k+1)^{-\alpha_{1}}
\sum_{l\in\mathbb{Z}_+}|w|^{2l}
(l+1)^{-\alpha_{2}}}\biggr)^{1/2}. \\
\end{align*}

Terms of the form
\[(k+1)^{-\beta_1}\]
are comparable to terms of the form
\[k^{-\beta_1},\]
which, for $k$ large, are comparable \cite[Chapter
5]{hardy:divergent} to terms of the form
\[\begin{pmatrix}
k-\beta_1\\
-\beta_1\\
\end{pmatrix}.
\]
The sum
\[\sum_{k\in\mathbb{Z}_+}\begin{pmatrix}
k-\beta_1\\
-\beta_1\\
\end{pmatrix}
|z|^{2k},\] with $1>\beta_1$, is comparable \cite[Chapter 5]{hardy:divergent} to
\[(1-|z|^2)^{\beta_1-1}.\]
We conclude that $|h(z,\,w)|$ is comparable to
\[\biggl((1-|z|^2)^{\frac{\beta_1-\alpha_1}{2}}(1-|w|^2)^{\frac{\beta_2-\alpha_2}{2}}\biggr).\\\]

Assume conversely that $|h(z,\,w)|$ is comparable to
\[\biggl((1-|z|^2)^{\frac{\beta_1-\alpha_1}{2}}(1-|w|^2)^{\frac{\beta_2-\alpha_2}{2}}\biggr),\]
and let $f$ be a function in $D_{\alpha}$. We recall that the
$D_{\beta}$ norm of $hf$ is equivalent to

\begin{align*}
\frac{1}{\pi^2}&\int_{0}^{1}\int_{0}^{1}\int_{0}^{2\pi}\int_{0}^{2\pi}
       |h(r_{1}e^{i\theta_{1}},\,r_{2}e^{i\theta_{2}})f(r_{1}e^{i\theta_{1}},\,r_{2}e^{i\theta_{2}})|^2\\
   &\hspace{1.5in} (1-r_{1}^2)^{-1-\beta_{1}}(1-r_{2}^2)^{-1-\beta_{2}}r_{1}r_{2}\,d\theta_{1}d\theta_{2}dr_{1}dr_{2}\\
   &\leq C\frac{1}{\pi^2}\int_{0}^{1}\int_{0}^{1}\int_{0}^{2\pi}
     \int_{0}^{2\pi}|f(r_{1}e^{i\theta_{1}},\,r_{2}e^{i\theta_{2}})|^2
     (1-r_1^2)^{\beta_1-\alpha_1}(1-r_2^2)^{\beta_2-\alpha_2}\\
   &\hspace{1.5in}(1-r_{1}^2)^{-1-\beta_{1}}(1-r_{2}^2)^{-1-\beta_{2}}r_{1}r_{2}\,d\theta_{1}d\theta_{2}dr_{1}dr_{2} \\
   &=\frac{C}{\pi^2}\int_{0}^{1}\int_{0}^{1}\int_{0}^{2\pi}\int_{0}^{2\pi}
       |f(r_{1}e^{i\theta_{1}},\,r_{2}e^{i\theta_{2}})|^2\\
   &\hspace{1.5in} (1-r_{1}^2)^{-1-\alpha_{1}}(1-r_{2}^2)^{-1-\alpha_{2}}r_{1}r_{2}\,d\theta_{1}d\theta_{2}dr_{1}dr_{2}.\\
\end{align*}
This last integral is  finite, as it is equivalent to the
$D_{\alpha}$ norm of $f$. We conclude that $h$ is in
$M(D_{\alpha},D_{\beta})$.
\end{proof}

\noindent If $\alpha\preceq(0,\,0)$ we have a complete
characterization of $\mathcal{M}(D_{\alpha})$.

\begin{prop}\label{prop:multipliers less than 0}
If $\alpha\preceq(0,\,0)$ then
$\mathcal{M}(D_{\alpha})=H^{\infty}(\mathbb{U}^2)$.
\end{prop}

\begin{proof}
We have seen that $\mathcal{M}(D_{\alpha})\subseteq
H^{\infty}(\mathbb{U}^2)$. It is clear from the integral
representation of the norm of $D_{\alpha}$, with
$\alpha\preceq(0,\,0)$, that we have
$H^{\infty}(\mathbb{U}^2)\subseteq\mathcal{M}(D_{\alpha})$.
\end{proof}

\noindent We also have a complete characterization of multipliers
from $D_{\alpha}$, $\alpha\succ (1,\,1)$, to $D_{\beta}$,
$\beta\preceq\alpha$.

\begin{thm}\label{thm:M D beta is D beta}
Let $\alpha\succ (1,\,1)$, $\beta\preceq\alpha$. Then
\[\mathcal{M}(D_{\alpha},\,D_{\beta})=D_{\beta}.\]
\end{thm}

\begin{proof}
We have the inclusion
$\mathcal{M}(D_{\alpha},\,D_{\beta})\subseteq D_{\beta}$. To prove
the inclusion
$D_{\beta}\subseteq\mathcal{M}(D_{\alpha},\,D_{\beta})$, we let
$f$ be an element of $D_{\alpha}$ and $g$ be an element of
$D_{\beta}$:
\begin{align*}
f(z,\,w)&=\sum_{(k,\,l)\in\mathbb{Z}^2_+}a_{k,\,l}z^kw^l, \textrm{ and }\\
g(z,\,w)&=\sum_{(k,\,l)\in\mathbb{Z}^2_+}b_{k,\,l}z^kw^l.
\end{align*}
Then
\[(f\cdot g)(z,\,w)=\sum_{(k,\,l)\in\mathbb{Z}^2_+}\biggl(\sum_{0\leq m\leq k}\sum_{0\leq n\leq l}
a_{m,\,n}b_{k-m,\,l-n}\biggr)z^kw^l,\] and
\[\|fg\|_{\beta}^2=\sum_{(k,\,l)\in\mathbb{Z}^2_+}(k+1)^{\beta_1}(l+1)^{\beta_2}\biggl|\sum_{0\leq m\leq k}\sum_{0\leq n\leq l}
a_{m,\,n}b_{k-m,\,l-n}\biggr|^2.\]

For ease of notation, let $\sum_{k,\,l}'$ denote $\sum_{0\leq
m\leq k}\sum_{0\leq n\leq l}$. We multiply each term,
$a_{m,\,n}b_{k-m,\,l-n}$, by
\[\frac{(m+1)^{\frac{\alpha_1}{2}}(n+1)^{\frac{\alpha_2}{2}}(k-m+1)^{\frac{\beta_1}{2}}(l-n+1)^{\frac{\beta_2}{2}}}
{(m+1)^{\frac{\alpha_1}{2}}(n+1)^{\frac{\alpha_2}{2}}(k-m+1)^{\frac{\beta_1}{2}}(l-n+1)^{\frac{\beta_2}{2}}},\]
and applying Cauchy-Schwarz inequality obtain
\begin{equation}\label{inequality:norm estimate}
\begin{split}
\|f&g\|_{\beta}^2\leq\sum_{(k,\,l)\in\mathbb{Z}^2_+}\biggl(\sideset{}{'}\sum_{k,\,l}
\frac{1}{(m+1)^{\alpha_1}(n+1)^{\alpha_2}(k-m+1)^{\beta_1}(l-n+1)^{\beta_2}}\biggr)\\
&\cdot \biggl(\sideset{}{'}\sum_{k,\,l}
(m+1)^{\alpha_1}(n+1)^{\alpha_2}(k-m+1)^{\beta_1}(l-n+1)^{\beta_2}
|a_{m,\,n}|^2|b_{k-m,\,l-n}|^2\biggr)\\
&\hspace{.5in}\cdot (k+1)^{\beta_1}(l+1)^{\beta_2}.
\end{split}
\end{equation}

Next, we notice that
\begin{equation}\label{inequality:k terms}
\begin{split}
(k+1)^{\beta_1}&\sum_{0\leq m\leq
k}\frac{1}{(m+1)^{\alpha_1}(k-m+1)^{\beta_1}}\\
&=\frac{(k+1)^{\beta_1}}{(k+2)^{\alpha_1}} \sum_{0\leq m\leq
k}\biggl(\frac{1}{m+1}+\frac{1}{k-m+1}\biggr)^{\alpha_1}\frac{1}{(k-m+1)^{\beta_1-\alpha_1}}\\
& \leq \sum_{0\leq m\leq
k}\biggl(\frac{1}{m+1}+\frac{1}{k-m+1}\biggr)^{\alpha_1}\\
 &\leq\biggl[\biggl(\sum_{0\leq m\leq
k}\frac{1}{(m+1)^{\alpha_1}}\biggr)^{\frac{1}{\alpha_1}}+\biggl(\sum_{0\leq
m\leq
k}\frac{1}{(k-m+1)^{\alpha_1}}\biggr)^{\frac{1}{\alpha_1}}\biggr]^{\alpha_1}\\
&\leq C_{\alpha_1},
\end{split}
\end{equation}
with $C_{\alpha_1}$ a constant depending only on $\alpha_1$. The
second last inequality follows from an application of Minkowski's
inequality. That the sums
\[\sum_{0\leq m\leq
k}\frac{1}{(m+1)^{\alpha_1}}\hspace{.2in}\textrm{ and
}\hspace{.2in}\sum_{0\leq m\leq k}\frac{1}{(k-m+1)^{\alpha_1}}\]
are bounded follows from the fact that $\alpha_1>1$.

Similarly, we obtain

\begin{equation}\label{inequality:l terms}
(l+1)^{\beta_2}\sum_{0\leq n\leq
l}\frac{1}{(n+1)^{\alpha_1}(l-n+1)^{\beta_2}}\leq C_{\alpha_2},
\end{equation}
with $C_{\alpha_2}$ a constant depending only on $\alpha_2$.

Since $m$ and $n$ are varying independently in the sum
\[(k+1)^{\beta_1}(l+1)^{\beta_2}\sideset{}{'}\sum_{k,\,l}
\frac{1}{(m+1)^{\alpha_1}(n+1)^{\alpha_2}(k-m+1)^{\beta_1}(l-n+1)^{\beta_2}},\]
we use Inequalities \ref{inequality:k terms} and \ref{inequality:l
terms} to estimate this product, and conclude that it is less than
$C_{\alpha_1}C_{\alpha_2}$.

Returning to Inequality \ref{inequality:norm estimate}, we see
that
\begin{align*}
\|fg\|_{\beta}^2&\leq
C_{\alpha_1}C_{\alpha_2}\sum_{(k,\,l)\in\mathbb{Z}^2_+}
\sideset{}{'}\sum_{k,\,l}
(m+1)^{\alpha_1}(n+1)^{\alpha_2}(k-m+1)^{\beta_1}(l-n+1)^{\beta_2}\\
&\hspace{1in}\cdot|a_{m,\,n}|^2|b_{k-m,\,l-n}|^2\\
&=C_{\alpha_1}C_{\alpha_2}\sum_{(k,\,l)\in\mathbb{Z}^2_+}
\sum_{(m,\,n)\in\mathbb{Z}^2_+}
(m+1)^{\alpha_1}(n+1)^{\alpha_2}(k+1)^{\beta_1}(l+1)^{\beta_2}\\
&\hspace{1in}\cdot|a_{m,\,n}|^2|b_{k,\,l}|^2\\
 &=C_{\alpha_1}C_{\alpha_2}\|f\|_{\alpha}^2\|g\|_{\beta}^2<\infty.
\end{align*}
Since $\|fg\|_{\beta}$ is finite, we have that $g$ is an element
of $\mathcal{M}(D_{\alpha},\,D_{\beta})$.
\end{proof}

\begin{cor}
If $\alpha\succeq\beta\succ(1,\,1)$, with $\alpha_1>\beta_1$ or $\alpha_2>\beta_2$,
then $\mathcal{M}(D_{\alpha})\subset\mathcal{M}(D_{\beta})$.
\end{cor}

\begin{proof}
This follows immediately from Theorem \ref{thm:M D beta is D beta}, Proposition
\ref{prop:proper containment}, and Remark \ref{rem:proper containment}.
\end{proof}

\noindent We have an analogous result, without proper inclusion,
for any $\alpha\succ\beta$.

\begin{thm}
If $\alpha \succ \beta$ then $M(D_{\alpha}) \subseteq
M(D_{\beta})$.
\end{thm}

\begin{proof}
If $f$ is in $M(D_{\alpha})$ then, by Corollary \ref{corollary:H
infty multipliers}, $f$ is in $H^\infty(\mathbb{U}^2)$.  By
Proposition \ref{prop:multipliers less than 0} we have that $f$ is
in $M(D_{\gamma})$ for all $\gamma \prec (0,0)$.  There exists a
$\lambda \in (0,1)$ such that
\[\beta - \lambda\alpha \prec(0,0).\]
Let $\gamma=(\gamma_1,\,\gamma_2)$ be defined as
\[\gamma=\frac{\beta -\lambda\alpha}{1-\lambda}.\]
Since $\gamma \prec (0,0)$ we have that $f$ is in $M(D_{\gamma})$.
Applying Theorem \ref{theorem:interpolation} with $\alpha$,
$\gamma$ and $\lambda$, we see that $f$ is in $D_{\beta}$.
\end{proof}

\noindent Thus far we have examined $\mathcal{M}(D_{\alpha},\,D_{\beta})$ with
$\alpha\succeq\beta$. The next theorem indicates why we have chosen to examine only these indices.

\begin{thm}\label{thm:multipliers alpha < beta}
If $\beta\succ\alpha$,
then $\mathcal{M}(D_{\alpha},\,D_{\beta})=\{0\}$.
\end{thm}

\begin{proof}
We begin by noting that
\[\|\lambda_{(z,\,w)}^{(\gamma_1,\,\gamma_2)}\|=\|\lambda_z^{\gamma_1}\|\|\lambda_w^{\gamma_2}\|.\]
The number on the left hand side is the norm of the point evaluation at $(z,\,w)\in\mathbb{U}^2$ on the
space $D_{(\gamma_1,\,\gamma_2)}$. The two numbers on the right hand side are the norm of the point evaluation
at $z\in\mathbb{U}$ (resp. $w\in\mathbb{U}$) on the space $D_{\gamma_1}$ (resp. $D_{\gamma_2}$).

Thus, by Lemma \ref{lemma:multiplier bound}, we have that if
$h\in\mathcal{M}(D_{\alpha},\,D_{\beta})$ then
\begin{align*}
|h(z,\,w)|&\leq \|T_h\|_{\alpha,\,\beta}\frac{\|\lambda_{(z,\,w)}^{\beta}\|}{\|\lambda_{(z,\,w)}^{\alpha}\|}\\
          &=\|T_h\|_{\alpha,\,\beta}\frac{\|\lambda_z^{\beta_1}\|\|\lambda_w^{\beta_2}\|}
{\|\lambda_{z}^{\alpha_1}\|\|\lambda_w^{\alpha_2}\|}.
\end{align*}

In the case of Dirichlet type spaces on the unit disc, Taylor
\cite[Theorem 4]{taylor:multiplier} shows that $a<b$ implies
$\mathcal{M} (D_{a},\,D_{b})=\{0\}$. An application of Taylor's
arguments shows that $\|\lambda_z^{\beta_1}\|/\|
\lambda_z^{\alpha_1}\|\rightarrow 0$ as $|z|\rightarrow 1$, and
that $\|\lambda_w^{\beta_2}\|/\| \lambda_w^{\alpha_2}\|\rightarrow
0$ as $|w|\rightarrow 1$. We conclude that $|h(z,\,w)|\rightarrow
0$ as $(z,\,w)\rightarrow \partial\mathbb{U}^2$. By the maximum
principle we conclude that $h\equiv 0$.

\end{proof}

\begin{rem}
It is of interest to note that, unlike most of the theorems we
prove, Theorem \ref{thm:multipliers alpha < beta} does not rely on
our ``partition'' of $\mathbb{R}^2$. That is to say, there is no
requirement that $\alpha_1$ and $\alpha_2$ lie in the same element
of Taylor's partition of $\mathbb{R}$. They are allowed to vary
independently of one another.

We note as well that a more general result than Theorem
\ref{thm:multipliers alpha < beta} can be proved (Proposition
\ref{prop:beta > alpha}), once we have examined slices of
functions in Section \ref{section:slices}.
\end{rem}

\section{Slices}\label{section:slices}
\noindent Given a function, $f$, on the bidisc, we can examine the
functions obtained by restricting $f$ to the one complex
dimensional slices of the bidisc parallel to the coordinate axes.
Specifically, we define $f_z$ and $f_w$ as follows.

\begin{definition}[$f_z$, $f_w$]
Let $f(z,\,w)$ be a function from $\mathbb{U}^2$ to $\mathbb{C}$,
and fix $w\in\mathbb{U}$. The function
$f_w(z):\mathbb{U}\rightarrow\mathbb{C}$ is defined as
\[f_w(z)=f(z,\,w).\]
Similarly, for fixed $z\in\mathbb{U}$ we define
$f_z(w):\mathbb{U}\rightarrow\mathbb{C}$ by
\[f_z(w)=f(z,\,w).\]
\end{definition}

\noindent It is natural to ask whether $f\in
D_{(\alpha_1,\,\alpha_2)}$ implies that $f_z\in D_{\alpha_2}$ and
$f_w\in D_{\alpha_1}$, for all $z,\,w\in\mathbb{U}$. We shall see
that the answer is yes, but that the converse is not necessarily
true.

\begin{thm}\label{thm:slices}
Let $f$ be an element of $D_{(\alpha_1,\,\alpha_2)}$. Then
$f_{w_0}\in D_{\alpha_1}$ for each $w_0\in\mathbb{U}$, and
$f_{z_0}\in D_{\alpha_2}$ for each $z_0\in\mathbb{U}$.
\end{thm}

\begin{proof}
We prove that $f_{w_0}$ is in $D_{\alpha_1}$. The  proof that
$f_{z_0}$ is in $D_{\alpha_2}$ is the same, with the obvious
modifications.

Fix $w_0 \in \mathbb{U}$. Let
\[f(z,\,w)=\sum_{(k,\,l)\in\mathbb{Z}^2_+}a_{k,\,l}z^kw^l\]
be an element of $D_{\alpha}$. We know that
\[\|f\|_{\alpha}^2=\sum_{(k,\,l)\in\mathbb{Z}^2_+}|a_{k,\,l}|^2(k+1)^{\alpha_1}(l+1)^{\alpha_2}\]
is finite. The function $f_{w_0}(z)$ is a holomorphic function
with power series expansion
\[f_{w_0}(z)=\sum_{k\in\mathbb{Z}_+}\biggl(\sum_{l\in\mathbb{Z}_+}a_{k,\,l}w_0^l\biggr)z^k.\]
For each $k\in\mathbb{Z}_+$ the  sum
$\sum_{l\in\mathbb{Z}_+}a_{k,\,l}w_0^l$ converges absolutely.

To show that $f_{w_0}$ is an element of $D_{\alpha_1}$, we must show that
\[\|f_{w_0}\|_{\alpha_1}^2=\sum_{k \in\mathbb{Z}_+}\biggl|\sum_{l\in\mathbb{Z}_+}a_{k,\,l}w_0^{l}\biggr|^2(k+1)^{\alpha_{1}} < \infty.\]
We have the following.
\begin{align*}
\sum_{k \in\mathbb{Z}_+}\biggl|\sum_{l \in\mathbb{Z}_+}a_{k,\,l}w_0^{l}&\biggr|^2(k+1)^{\alpha_{1}} \leq
\sum_{k \in\mathbb{Z}_+}\biggl(\sum_{l \in\mathbb{Z}_+}|a_{k,\,l}||w_0|^{l}\biggr)^2(k+1)^{\alpha_{1}} \\
&\leq (1-|w_0|)^{-1}\sum_{k \in\mathbb{Z}_+}\biggl(\sum_{l \in\mathbb{Z}_+}|a_{k,\,l}|^2|w_0|^{l}\biggr)(k+1)^{\alpha_{1}} \\
&\leq (1-|w_0|)^{-1}\biggl(C_{w_0} + \sum_{k \in\mathbb{Z}_+}\sum_{l \in\mathbb{Z}_+}|a_{k,\,l}|^2(k+1)^{\alpha_{1}}(l+1)^{\alpha_{2}}\biggr)\\
&=(1-|w_0|)^{-1}(C_{w_0}+\|f\|_{\alpha}^2)<\infty.
\end{align*}
The second inequality is an application of Jensen's Inequality,
and the fact that $\sum_{l\in\mathbb{Z}_+}|w_0|^l=(1-|w_0|)^{-1}$.
The last inequality follows from reasoning similar to that
employed in the proof of Lemma \ref{lemma:multiindex convergence}.
The constant $C_{w_0}$ is a finite number depending on $w_0$.
\end{proof}

\noindent As mentioned above, the converse to Theorem
\ref{thm:slices} is not true.

\begin{thm}\label{thm:non factoring}
Fix $\alpha=(\alpha_{1},\alpha_{2})$.  Then there exists a
function $f:\mathbb{U}^2 \rightarrow \mathbb{C}$ such that $f_{w}
\in D_{\alpha_{1}}$ for all $w$ in $\mathbb{U}$, $f_{z} \in
D_{\alpha_{2}}$ for all $z$ in $\mathbb{U}$, but $f$ is not in
$D_{\alpha}$.
\end{thm}

\begin{proof}
Let $f:\mathbb{U}^2 \rightarrow \mathbb{C}$ be defined as
$f(z,\,w)=\sum_{(k,\,l)\in\mathbb{Z}^2_+}a_{k,\,l}z^l w^k$, where
\[a_{k,\,l} = \sqrt{\frac{(k+1)^{1-\alpha_{1}}(l+1)^{1-\alpha_{2}}}{(k+1)^{3}+ (l+1)^3}}.\]

We fix $w_0$ in $\mathbb{U}$. Then
\begin{align*}
\|f_{w_0}\|^2_{\alpha_{1}} &= \sum_{k\in\mathbb{Z}_+} \biggl|
\sum_{l \in\mathbb{Z}_+}
\biggl(\frac{(k+1)^{1-\alpha_{1}}(l+1)^{1-\alpha_{2}}}{(k+1)^{3}+ (l+1)^3}\biggr)^{\frac{1}{2}}w_0^l\biggr|^2(k+1)^{\alpha_{1}} \\
&\leq
(1-|w_0|)^{-1}\sum_{(k,\,l)\in\mathbb{Z}^2_+}\frac{(k+1)(l+1)}{(k+1)^{3}+
(l+1)^3}(l+1)^{-\alpha_{2}}|w_0|^l
\end{align*}
by an application of Jensen's inequality.

Note that there is an $L$ such that for each $k\in\mathbb{Z}_+$,
and for $l\geq L$, we have that
\begin{align*}
|w_0|^l &< (l+1)^{\alpha_{2}-3} + (k+1)^{-3}(l+1)^{\alpha_{2}} \\
&= (k+1)^{-2}(l+1)^{\alpha_{2}-2}\biggl(\frac{(k+1)^{3}+
(l+1)^3}{(k+1)(l+1)}\biggr).
\end{align*}
Thus, for each $k\in\mathbb{Z}_+$, and for $l\geq L$,
\[\frac{(k+1)(l+1)}{(k+1)^{3}+ (l+1)^3}(l+1)^{-\alpha_{2}}|w_0|^l < (k+1)^{-2}(l+1)^{-2}.\]

We conclude that
 \begin{align*}
\|f_{w_0}\|^2_{\alpha_{1}} &\leq (1-|w_0|)^{-1}
\biggl(\sum_{k \in\mathbb{Z}_+}\sum_{l\geq L}(k+1)^{-2}(l+1)^{-2}\\
&+ \sum_{0\leq l\leq L}\sum_{k
\in\mathbb{Z}_+}\frac{(k+1)(l+1)}{(k+1)^{3}+
(l+1)^3}(l+1)^{-\alpha_{2}}|w_0|^l \biggr).
\end{align*}
The first summand is clearly finite. The second summand is finite,
since for each $l\in\{0,\,\ldots,\,L\}$ the series
\[\sum_{k \in\mathbb{Z}_+}\frac{(k+1)(l+1)}{(k+1)^{3}+ (l+1)^3}(l+1)^{-\alpha_{2}}|w_0|^l\]
is dominated by a constant multiple of
$\sum_{k\in\mathbb{Z}_+}\frac{1}{(k+1)^2}$. Hence
$\|f_{w_0}\|_{\alpha_1}$ is finite, and $f_{w_0}$ is in
$D_{\alpha_{1}}$.

We have shown that $f_w$ is in $D_{\alpha_1}$, for each
$w\in\mathbb{U}$. By the symmetry of the Taylor coefficients of
$f$  we see that $f_{z}$ is in $D_{\alpha_{2}}$, for each $z$ in
$\mathbb{U}$.

Since $f_z$ is in $D_{\alpha_2}$ for each $z\in\mathbb{U}$, we see
that $f_z$ is holomorphic for each $z\in\mathbb{U}$. Similarly,
$f_w$ is holomorphic for each $w\in\mathbb{U}$. We conclude that
$f$ is holomorphic on $\mathbb{U}^2$.

The $D_{\alpha}$ norm of $f$ is, however, infinite:
\[\|f\|_{\alpha}^2=\sum_{(k,\,l)\in\mathbb{Z}^2_+}\frac{(k+1)(l+1)}{(k+1)^{3}+ (l+1)^3} >
\sum_{k\in\mathbb{Z}_+}\frac{1}{2(k+1)} = \infty.\]

\end{proof}

\begin{rem} We remark that analogous examples can be found in
higher dimensions. Let $f:\mathbb{U}^n\rightarrow\mathbb{C}$ be
defined as
\[f(z_1,\,\ldots,\,z_n)=\sum_{j_1\in\mathbb{Z}_+}\cdots\sum_{j_n\in\mathbb{Z}_+}
a_{j_1,\,\ldots,\,j_n}z_1^{j_1}\cdots z_n^{j_n},\] with
\[a_{j_1,\,\ldots,\,j_n}=\biggl(\frac{j_1^{1-\alpha_1}\cdots
j_n^{1-\alpha_n}}{j_1^{2n-1}+\cdots+j_n^{2n-1}}\biggr)^{\frac{1}{2}}.\]
Then $f$ is not in $D_{(\alpha_1,\,\ldots,\,\alpha_n)}$. However,
if we fix $k$ coordinates, $1\leq k\leq n-1$, then $f$ is in the
corresponding $n-k$ dimensional Dirichlet space. For example, fix
$z_1=b_1,\,\ldots,\,z_k=b_k$, with $k\in\{1,\,\ldots,\,n-1\}$.
Then
\[f(b_1,\,\ldots,\,b_k,\,z_{k+1},\,\ldots,\,z_n):\mathbb{U}^{n-k}
\rightarrow\mathbb{C}\] is in
$D_{(\alpha_{k+1},\,\ldots,\,\alpha_n)}$.

\end{rem}

\noindent Theorem \ref{thm:slices} tells us that $f_w$ is in
$D_{\alpha_1}$ for each $w\in\mathbb{U}$. It does not, however,
indicate whether the norms $\|f_w\|_{\alpha_1}$ are uniformly
bounded in $w$.

An example illustrates that in general there need be no such
bound. Consider the function
$f(z,\,w)=\sum_{(k,\,l)\in\mathbb{Z}^2_+}a_{k,\,l}z^kw^l$, where
$a_{k,\,l}=1$ for all $(k,\,l)\in\mathbb{Z}^2_+$. Let
$\alpha_1,\,\alpha_2<-1$. The function $f$ is in
$D_{(\alpha_1,\,\alpha_2)}$, with
\[\|f\|_{(\alpha_1,\,\alpha_1)}^2=\sum_{(k,\,l)\in\mathbb{Z}^2}(k+1)^{\alpha_1}(l+1)^{\alpha_2}.\]
On the other hand, the $D_{\alpha_1}$ norm of $f_w$ is
\begin{align*}
\|f_w\|^2_{\alpha_1}&=\sum_{k\in\mathbb{Z}_+}\biggl|\sum_{l\in\mathbb{Z}_+}w^l\biggr|^2(k+1)^{\alpha_1}\\
              &=\biggl|\frac{1}{1-w}\biggr|^2\sum_{k\in\mathbb{Z}_+}(k+1)^{\alpha_1}\\
              &=\biggl|\frac{1}{1-w}\biggr|^2 C.
\end{align*}
Clearly $\|f_w\|_{\alpha_1}\rightarrow\infty$ as $w\rightarrow 1$.

This example motivates the following question regarding a partial
converse to Theorem \ref{thm:slices}. Let $f$ be a function,
$f:\mathbb{U}^2\rightarrow\mathbb{C}$. Assume that
$\|f_w\|_{\alpha_1}<M$ for each $w\in\mathbb{U}$. When does this
imply that $f\in D_{(\alpha_1,\,\alpha_2)}$? For which $\alpha_2$?

We present several results in this direction. The first result
indicates that, in general, uniform boundedness in $w$ of
$\|f_w\|_{\alpha_1}$ is not enough to guarantee that $f$ is in
$D_{(\alpha_1,\,\alpha_2)}$. The second result, however, indicates
that if we restrict our attention to a suitable set of
$(\alpha_1,\,\alpha_2)$, then uniform boundedness in $w$ of
$\|f_w\|_{\alpha_1}$ does guarantee that $f$ is in
$D_{(\alpha_1,\,\alpha_2)}$. Finally, we illustrate a simple
situation where a converse to Theorem \ref{thm:slices} holds.

\begin{thm}
Let $\alpha_2>1$. For any $\alpha_1\in\mathbb{R}$, there is a
function $f:\mathbb{U}^2\rightarrow\mathbb{C}$ such that
$\|f_w\|_{\alpha_1}<M<\infty$ for all $w\in\mathbb{U}$, but $f$ is
not in $D_{(\alpha_1,\,\alpha_2)}$.
\end{thm}

\begin{proof}
Since $\alpha_2>1$ we know that $D_{\alpha_2}\subseteq
H^{\infty}(\mathbb{U})$. Since $H^{\infty}(\mathbb{U})$ is not a
Hilbert space, while $D_{\alpha_2}$ is, the containment is strict.
Let $g$ be a function in $H^{\infty}(\mathbb{U})\backslash
D_{\alpha_2}$, and assume that $\|g\|_{H^{\infty}}=1$.

Define
\[g_k=\frac{g}{(k+1)^{\frac{\alpha_1+2}{2}}},\]
so that
\[\|g_k\|_{\infty}=\frac{1}{(k+1)^{\frac{\alpha_1+2}{2}}},\]
 and expand $g_k$ as
\[g_k(w)=\sum_{l\in\mathbb{Z}_+}a_{k,\,l}w^l.\]
Clearly $g_k$ is not in $D_{\alpha_2}$.

\noindent Define $f(z,\,w)$ as
\[f(z,\,w)=\sum_{k,\,l\in\mathbb{Z}^2_+}a_{k,\,l}z^kw^l.\]

\noindent Then
\begin{align*}
\|f\|^2_{(\alpha_1,\,\alpha_2)}&=\sum_{k,\,l\in\mathbb{Z}^2_+}|a_{k,\,l}|^2(k+1)^{\alpha_1}
(l+1)^{\alpha_2}\\
&=\sum_{k\in\mathbb{Z}_+}\biggl(\sum_{l\in\mathbb{Z}_+}|a_{k,\,l}|^2(l+1)^{\alpha_2}\biggr)(k+1)^{\alpha_1}\\
&=\sum_{k\in\mathbb{Z}_+}\|g_k\|^2_{\alpha_2}(k+1)^{\alpha_1}\\
 &=\infty.
\end{align*}

On the other hand,
\begin{align*}
\|f_w\|^2_{\alpha_1}&=\sum_{k\in\mathbb{Z}_+}\biggl|\sum_{l\in\mathbb{Z}_+}a_{k,\,l}
w^l\biggr|^2
(k+1)^{\alpha_1}\\
&\leq \sum_{k\in\mathbb{Z}_+}(k+1)^{-2}\\
&=C,
\end{align*}
where $C$ is a finite constant independent of $w$.
\end{proof}

\noindent While we cannot guarantee that uniform boundedness in $w$ of
$\|f_w\|_{\alpha_1}$ implies inclusion in
$D_{(\alpha_1,\,\alpha_2)}$, we now see that there are situations
where this is the case.

\begin{thm}
Let $\alpha_1\leq 0$. Suppose that
$f:\mathbb{U}^2\rightarrow\mathbb{C}$ is such that
$\|f_w\|_{\alpha_1}<M<\infty$ for all $w\in\mathbb{U}$. Then, for
any $\alpha_2\leq 0$, $f$ is in $D_{(\alpha_1,\,\alpha_2)}$.
\end{thm}

\begin{proof}
Since $f_w$ is in $D_{\alpha_1}$ and $\alpha_1\leq 0$, by Lemma
\ref{lemma:integral norm equivalence} we have that
$\|f_w\|_{\alpha_1}$ is comparable to
\[I_{\alpha_1}(f_w)=\frac{1}{\pi}\int_{0}^{2\pi}\int_0^1|f_w(re^{i\theta})|^2(1-r)^{-1-\alpha_1}r\,drd\theta\]
if $\alpha_1<0$, and comparable to
\[I_{\alpha_1}(f_w)=\sup_{0\leq r
<1}\int_0^{2\pi}|f_w(re^{i\theta})|^2\,d\theta\] if $\alpha_1=0$.

We are assuming that $\|f_w\|_{\alpha_1}<M<\infty$ for all
$w\in\mathbb{U}$, so in fact $I_{\alpha_1}(f_w)<M$ for all
$w\in\mathbb{U}$.

For ease of notation we assume that $\alpha_1$, $\alpha_2$ are
both less than $0$. In this case $\|f\|_{(\alpha_1,\,\alpha_2)}$
is equivalent to
\begin{align*}
&\frac{1}{\pi}\int_{0}^{2\pi}\int_0^1\biggl(\frac{1}{\pi}\int_{0}^{2\pi}\int_0^1
|f(r_1e^{i\theta_1},\,r_2e^{i\theta_2})|^2(1-r_1)^{-1-\alpha_1}r_1\,dr_1d\theta_1\biggr)\\
&\hspace{2in}(1-r_2)^{-1-\alpha_2}r_2\,dr_2d\theta_2\\
&<\frac{1}{\pi}\int_{0}^{2\pi}\int_0^1M(1-r_2)^{-1-\alpha_2}r_2\,dr_2d\theta_2\\
&< \infty.
\end{align*}
We conclude that $f$ is in $D_{(\alpha_1,\,\alpha_2)}$.

(We remark that if $\alpha_1$ or $\alpha_2$ are equal to $0$, then
the obvious adjustments must be made, but the proof goes through
in exactly the same fashion.)
\end{proof}

\noindent While it is in general not the case that the converse to
Theorem \ref{thm:slices} is true, in the following simple
situation it does hold.

\begin{prop}
Let $f_{j}:\mathbb{U} \rightarrow \mathbb{C}$ be an element of
$D_{\alpha_{j}}$, for $j=1,\, 2$. Then $f(z,\,
w)=f_{1}(z)f_{2}(w)$ is in $D_{(\alpha_{1},\, \alpha_{2})}$.
\end{prop}

\begin{proof}
The Taylor coefficients of $f$ are
\[a_{k,\,l}=b^1_kb^2_l,\]
where $b^j_k$ denotes the $k^{th}$ Taylor coefficient of $f_{j}$.
We see that
\begin{align*}
\|f\|_{\alpha}^2&=\sum_{(k,\,l)\in\mathbb{Z}^2_+}|a_{k,\,l}|^2(k+1)^{\alpha_{1}}(l+1)^{\alpha_{2}}\\
&= \sum_{k \in\mathbb{Z}_+}|b^1_k|^2(k+1)^{\alpha_{1}}\sum_{l
\in\mathbb{Z}_+}|b^2_l|^2(l+1)^{\alpha_{2}}\\
&=\|f_1\|_{\alpha_1}^2\|f_2\|_{\alpha_2}^2 < \infty,
\end{align*}
and conclude that $f$ is in $D_{(\alpha_1,\,\alpha_2)}$.
\end{proof}

\noindent Having examined the restriction of functions in
$D_{\alpha}$ to  lower dimensional discs parallel to the
coordinate axes, we do the same with multipliers. Just as the
restriction, $f_w$, of a function was in the corresponding lower
dimensional Dirichlet type space, so the restriction of a
multiplier is a multiplier of the lower dimensional space.

\begin{thm}\label{thm:multiplier slices}
Let $h$ be an element of
$\mathcal{M}(D_{(\alpha_1,\,\alpha_2)},\,D_{(\beta_1,\,\beta_2)})$.
Then\\ $h_{w_0}\in \mathcal{M}(D_{\alpha_1},\,D_{\beta_1})$ for
each $w_0\in\mathbb{U}$, and $h_{z_0}\in
\mathcal{M}(D_{\alpha_2},\,D_{\beta_2}) $ for each
$z_0\in\mathbb{U}$.
\end{thm}

\begin{proof}
We prove that $h_{w_0}$ is in
$\mathcal{M}(D_{\alpha_1},\,D_{\beta_1})$. The statement regarding
$h_{z_0}$ is proved in the same fashion.

Let $h$ be an element of
$\mathcal{M}(D_{(\alpha_1,\,\alpha_2)},\,D_{(\beta_1,\,\beta_2)})$,
and fix $w_0\in\mathbb{U}$. Let
$f(z):\mathbb{U}\rightarrow\mathbb{C}$ be an  element of
$D_{\alpha_1}$, and define $\widetilde{f}(z,\,w):
\mathbb{U}^2\rightarrow\mathbb{C}$ by
\[\widetilde{f}(z,\,w)=f(z).\]
Notice that $\widetilde{f}$ is an element of
$D_{(\alpha_1,\,\alpha_2)}$; in fact, $\widetilde{f}$ is an
element of $D_{(\alpha_1,\,\beta)}$, for any $\beta\in\mathbb{R}$.

Since $h$ is in
$\mathcal{M}(D_{(\alpha_1,\,\alpha_2)},\,D_{(\beta_1,\,\beta_2)})$,
we have that $h\widetilde{f}$ is in $D_{(\beta_1,\,\beta_2)}$. We
note that $(h\widetilde{f})_{w_0}=
h_{w_0}\widetilde{f}_{w_0}=h_{w_0}f$, and by Theorem
\ref{thm:slices} we have that $(h\widetilde{f})_{w_0}$ is in
$D_{\beta_1}$. We conclude that $h_{w_0}$ is an element of
$\mathcal{M}(D_{\alpha_1},\,D_{\beta_1})$.
\end{proof}

\noindent Just as the converse to Theorem \ref{thm:slices} was not
in general true, we shall see that the converse to Theorem
\ref{thm:multiplier slices} is also not true in general. In fact,
the counterexample to the converse of Theorem \ref{thm:slices}
(Theorem \ref{thm:non factoring}) is the key ingredient in showing
that this is the case.

\begin{prop}
Let $\gamma\succ(1,\,1)$, and let $\alpha\preceq\gamma$. Then
there is a function, $f$, such that
$f_{w_0}\in\mathcal{M}(D_{\gamma_1},\,D_{\alpha_1})$, and
$f_{z_0}\in\mathcal{M}(D_{\gamma_2},\,D_{\alpha_2})$, for each
$z_0,\,w_0\in\mathbb{U}$, but $f$ is not in
$\mathcal{M}(D_{\gamma},\,D_{\alpha})$.
\end{prop}

\begin{proof}
By Theorem \ref{thm:M D beta is D beta} we have that
$\mathcal{M}(D_{\gamma},\,D_{\alpha})=D_{\alpha}$. Let $f$ be as
in the proof of Theorem \ref{thm:non factoring}, so that
$f_{w_0}\in D_{\alpha_1}$, $f_{z_0}\in D_{\alpha_2}$, for each
$z_0,\,w_0\in\mathbb{U}$, but $f$ is not in $D_{\alpha}$.

We know that $f_{w_0}\in\mathcal{M}(D_{\gamma_1},\,D_{\alpha_1})$
and $f_{z_0}\in\mathcal{M}(D_{\gamma_2},\,D_{\alpha_2})$, by
Taylor's \cite[Theorem 7]{taylor:multiplier} one variable  version
of Theorem \ref{thm:M D beta is D beta}. The function $f$,
however, is not in $\mathcal{M}(D_{\gamma},\,D_{\alpha})$, as it
is not an element of $D_{\alpha}$.
\end{proof}

\noindent We present one final example, showing further evidence
that a converse to Theorem \ref{thm:multiplier slices} does not in
general hold.

\begin{example}
Let $f(z,\,w)=\frac{1}{z+w-2}$. Clearly $f$ is not in
$H^{\infty}(\mathbb{U}^2)$, and thus $f$ is not in
$\mathcal{M}(D_{\alpha})$ for any $\alpha\in\mathbb{R}^2$. The
slice functions, $f_{w_0}$ and $f_{z_0}$, however, are in
$H^{\infty}(\mathbb{U})$ for each $z_0,\,w_0\in\mathbb{U}$. Thus,
if $\alpha\preceq(0,\,0)$, then
$f_{z_0}\in\mathcal{M}(D_{\alpha_2})$ and
$f_{w_0}\in\mathcal{M}(D_{\alpha_1})$.\\

We remark as well that $f$ is not in $D_{\alpha}$ for any
$\alpha\succeq(0,\,0)$. To see this we note that $f$ has the power
series expansion
\[f(z,\,w)=-\sum_{k\in\mathbb{Z}_+}\frac{1}{2^{k+1}}\sum_{l=0}^k
\begin{pmatrix}
k\\
l\\
\end{pmatrix}
z^lw^{k-l}.\]
Thus
\begin{align*}
\|f\|_{\alpha}^2&=\sum_{k\in\mathbb{Z}_+}\frac{1}{2^{2(k+1)}}\sum_{l=0}^k
 \begin{pmatrix}
 k\\
 l\\
 \end{pmatrix}^2(l+1)^{\alpha_1}(k-l+1)^{\alpha_2}\\
&\geq\sum_{k\in\mathbb{Z}_+}\frac{k+1}{2^{2(k+1)}}\biggl(\sum_{l=0}^k
\begin{pmatrix}
k\\
l\\
\end{pmatrix}^2\frac{1}{k+1}\biggr)\\
&\geq\sum_{k\in\mathbb{Z}_+}\frac{k+1}{2^{2(k+1)}}\biggl(\sum_{l=0}^k
\begin{pmatrix}
k\\
l\\
\end{pmatrix}\biggr)^2\frac{1}{(k+1)^2}\\
&=\sum_{k\in\mathbb{Z}_+}\frac{1}{4(k+1)}=\infty.
\end{align*}

We claim, however, that  $f_{w_0}$ is in $D_{\alpha_1}$, for any
$w_0\in\mathbb{U}$, and any $\alpha_1\in\mathbb{R}$. To see this
we first note that we can  expand $f(z,\,w)$ as
\[f_w(z)=-\sum_{k\in\mathbb{Z}_+}\biggl(\frac{1}{2-w}\biggr)^{k+1}z^k.\]
We notice that $1/(|2-w_0|)<1$. Thus there is an $N$ such that if
$n\geq N$ then
\[\biggl(\frac{1}{|2-w_0|}\biggr)^{2(n+1)}<(n+1)^{-\alpha_1-2}.\]
Then
\begin{align*}
\|f_{w_0}\|_{\alpha_1}^2&=\sum_{k\in\mathbb{Z}_+}\biggl(\frac{1}{|2-w_0|}\biggr)^{2(k+1)}(k+1)^{\alpha_1}\\
                      &\leq\sum_{k=0}^{N}\biggl(\frac{1}{|2-w_0|}\biggr)^{2(k+1)}(k+1)^{\alpha_1}
                      +\sum_{k=N+1}^{\infty}(n+1)^{-2}<\infty,\\
\end{align*}
and thus $f_{w_0}$ is in $D_{\alpha_1}$.

 Analogous calculations show that $f_{z_0}$ is in
$D_{\alpha_2}$, for any $\alpha_2\in\mathbb{R}$, and any
$z_0\in\mathbb{U}$.

Fix $\gamma\succ(1,\,1)$, and
$(0,\,0)\preceq(\alpha_1,\,\alpha_2)\preceq(\gamma_1,\,\gamma_2)$.
We saw above that $f$ is not in $D_{\alpha}$, and is thus not in
$\mathcal{M}(D_{\gamma},\,D_{\alpha})$. We have also seen that
$f_{z_0}$ is in $D_{\alpha_2}$, and is thus in
$\mathcal{M}(D_{\gamma_2},\,D_{\alpha_2})$, and that $f_{w_0}$ is
in $D_{\alpha_1}$, and is thus in
$\mathcal{M}(D_{\gamma_1},\,D_{\alpha_1})$.
\end{example}

\noindent We mention, finally, a simple corollary of our work with
slices. We easily generalize Theorem \ref{thm:multipliers alpha <
beta}.

\begin{prop}\label{prop:beta > alpha}
Let $\beta$ and $\alpha$ be elements of $\mathbb{R}^2$. If
$\beta_1>\alpha_1$ then
$\mathcal{M}(D_{\alpha},\,D_{\beta})=\{0\}$, regardless of the
relationship between $\alpha_2$ and $\beta_2$.

The same is true if $\beta_2>\alpha_2$, regardless of the
relationship between $\alpha_1$ and $\beta_1$.
\end{prop}

\begin{proof}
We recall that Taylor \cite[Theorem 4]{taylor:multiplier} proved
that if $a<b$ then $\mathcal{M}(D_a,\,D_b)=\{0\}$.

Assume that $\alpha_1<\beta_1$. We know from Theorem
\ref{thm:multiplier slices} that if $f$ is in
$\mathcal{M}(D_{\alpha},\,D_{\beta})$ then $f_w(z)$ is in
$\mathcal{M}(D_{\alpha_1},\,D_{\beta_1})$ for each
$w\in\mathbb{U}$. This implies that $f_w$ is identically zero for
each $w\in\mathbb{U}$, and thus that $f\equiv0$.
\end{proof}

\subsection*{Acknowledgements}
The authors would like to thank Professor Ron Douglas for providing
valuable suggestions during the writing of this paper.  The authors
would also like to thank the referee for his comments.

\end{document}